\documentclass[11pt]{article}
\usepackage{myStyleSheet}

\begin{document}

\title{Hypergraphic Oriented Matroid Relational Dependency Flow Models of Chemical Reaction Networks}
\author{Colin G. Bailey, Dean W. Gull\footnote{corresponding author: gulld@amath.washington.edu}, and Joseph S. Oliveira}
\bibliographystyle{plain}
\maketitle


\begin{abstract}
In this paper we derive and present an application of hypergraphic oriented matroids for the purpose of enumerating the variable interdependencies that define the chemical complexes associated with the kinetics of non-linear dynamical system representations of chemical kinetic reaction flow networks.  The derivation of a hypergraphic oriented matroid is obtained by defining a closure operator on families of $n$-subsets of signed multi-sets from which a ``$\mathbb{Z}$-module'' is obtained.  It has been observed that every instantiation of the closure operator on the signed multiset families define a matroid structure.  It is then demonstrated that these structures generate a pair of dual matroids corresponding respectively to hyperspanning trees and hypercycles obtained from the corresponding directed hypergraphs.  These structures are next systematically evaluated to obtain solution sets that satisfy systems of non-linear chemical kinetic reaction flow networks in the MAP Kinase cascade cell-signaling pathway.
\end{abstract}

\section{Introduction}
Numerous algebraic-geometric combinatorial models have been developed to analyze chemical reaction systems.  Sellers \cite{sellers84} developed an algebraic-combinatorial model for identifying steady state and ``cycle-free" solutions to discrete representations of derived dynamic systems, in which the reaction dynamics are given in terms of a linear transformation that defines set intersection families of hyperplanes.  The signed half-spaces implied by families of sets of hyperplane arrangements, in turn, generate a signed convex polyhedron that characterizes all possible flows modulo scalar multiples.  Goss and Peccoud \cite{goss98} demonstrated the use of Petri nets for studying chemical reaction systems.  Oliveira \emph{et al.} \cite{oliveira01} showed that the Petri nets used by Goss and Peccoud are equivalent to directed hypergraphs (hyperdigraphs) and proposed using oriented matroids to analyze the topological structure of the hyperdigraph derived from a given chemical reaction system.  Further, Oliveira \emph{et al.} \cite{oliveira03,oliveira04} demonstrated the utility of this topological hypergraph analysis using the Krebs cycles, and then later in a sequel the EGFR signaling network.  Yang \emph{et al.} \cite{yang05} also demonstrated the use of oriented matroids in analyzing steady state solutions to chemical reaction dynamics.  

The theory of oriented matroids provides a combinatorial algebraic-geometric abstraction of linear spaces based on signed orientations \cite{aigner97,bjorner99}.  In all generality, the stoichiometric number, which corresponds to the molecular quantity of chemical reactant required to drive a reaction, will be greater than or equal to zero, and not restricted to the set $\{-1,0,+1\}$.  Therefore, we consider the stoichiometric number as being a ``ranking'' edge weight.  This is the reason we have elected to define the oriented matroid representation of the chemical reaction system over the $\mathbb{Z}$-module which maintains the multiplicities.

Multisets are used to track the multiplicity of member elements in a given base set $S$.  Based on this notion, signed multisets are instantiated multisets with the added generalization that elements of the base set are allowed to have ``negative multiplicity".  Multisets and signed multisets are used extensively in the study of hypergraphs and Petri nets \cite{blizard89,reisig85} as a representation for vertex markings and edge weights.  In the latter case, \cite{reisig85} uses the term ``multirelation", when the multiset base set is given as a relation $S \subseteq D\times D$ for some set $D$.  We demonstrate that multisets are generalizable as a $\mathbb{Z}$-module, over which we will construct a hyperdigraph model of chemical reaction networks.    An observable consequence of the corresponding matroid structure is a constructive method for obtaining ``hyper-spanning" forests and enumerating the hypercycles in hyperdigraphs with weights \cite{calude01,reisig85}.  This structure will enable us to qualitatively analysis of the steady state conditions in the chemical reaction network.

Every hypergraph is a generalization of a graph \cite{berge73} in which each hyperedge may be incident with more than two vertices.  The motivation for using hypergraphs as a representation of the chemical reaction network comes from the fact that not all of the reactions in such systems necessarily need to be unimolecular.  A unimolecular reaction involves a single chemical species undergoing a molecular change to produce a different chemical species.  In general, many chemical species react as a chemical complex to produce a different chemical complex, which is in turn composed of many chemical species.  We will therefore model the chemical reactions as directed hyperedges, thus generating a hyperdigraph as a model of the chemical reaction network.

We have validated and verified this discrete systems model approach by introducing the analysis of the mitogen-activated protein kinase cascade (MAP Kinase cascade) to demonstrate the effective use of applying the hyperdigraph model of chemical reaction networks \cite{bhalla99}.  The MAP Kinase cascade signaling pathway is a phosphorylation cascade that activates a set of cell regulatory molecules that govern cellular processes such as cell proliferation, differentiation, and development.  It is interesting and worth noting from an evolutionary perspective that many distinct and disjoint cellular signaling processes have been identified that utilize the MAP Kinase cascade pattern.  The MAP Kinase cascade is therefore considered to be an isoform for this category of cascading biochemical reaction networks \cite{seger95}.  The hyperdigraph representation provides a novel approach for the analysis of the steady state flux conditions.  Further, we will utilize the hyperdigraph model of the MAP Kinase cascade to address an open question by \cite{estrada06} regarding the graphic centrality rank index analysis of hyperdigraph models of chemical reaction systems.

\section{Hyperdigraph Model of Chemical Reaction Networks}
First let $S$ be a finite set of chemical species over which we will define chemical complexes.  A chemical complex $c \in \mathbb{N}^S$ is a multiset defined on the set of chemical species, where the multiplicity $c(s)$ $(s \in S)$ is the molecularity of species $s$ in the given complex.  We denote this finite collection of chemical complexes by $C \subseteq \mathbb{N}^S$.  The finite set of reactions then is defined by the relation $\mathcal{R} \subseteq C \times C$.  Next define two projection functions $\rho$ and $\pi$ from $C \times C$ into $C$, where $\rho(c,c') \mapsto c$ is to be referred to as the reactant and correspondingly, $\pi(c,c') \mapsto c'$ is the product of the reaction defined by the relation $c \mathcal{R} c'$.

Now we define two matrices $A,B \in \mathbb{N}^{\mathcal{R} \times S}$, where $r \equiv (c,c') \in \mathcal{R}$ and $s \in S$, such that
\[
A(r, s) = \rho(c,c')(s) = c(s),
\]
and
\[
B(r, s) = \pi(c,c')(s) = c'(s).
\]
Observe now that the matrices $A$ and $B$ correspond to the molecularities of the reactant and product chemical complexes, respectively, for each reaction.  The net change in the species molecularity, of the given complex, can be modeled by the matrix $N \in \mathbb{Z}^{S\times \mathcal{R}}$ where
\[
N = (B - A)^T.
\]
Further observe that for $s \in S$ and $r \in \mathcal{R}$ it follows that
\[
N(s,r) = B^T(s,r) - A^T(s,r) = \pi(r)(s) - \rho(r)(s).
\]
Therefore the columns of the matrix $N$ are presented as a finite collection of signed multisets defined on the set $S$.

A hyperdigraph $\mathcal{H} = (V,\mathcal{E})$ is defined by a finite vertex set $V$ and an edge set $\mathcal{E} \subseteq \{-1,0,+1\}^V$ such that each hyperedge $E \in \mathcal{E}$ is a signed set on the vertex set $V$ so that $E$ partitions the set $V$ into classes $E^+, ~E^-$, and $E^0$.  A weighted hyperdigraph $\mathcal{H} = (V,\mathcal{E},W)$ is a hyperdigraph $\mathcal{H} = (V,\mathcal{E})$ with an assigned weight function given by $W \in \mathbb{R}_+^{V\times \mathcal{E}}$.  Evidently the incidence matrix of the hyperdigraph $M \in \mathbb{R}^{V \times \mathcal{E}}$ is defined for $v \in V$ and $E \in \mathcal{E}$ by setting
\[
M(v, E) = E(v)W(v,E) = \left\{\begin{array}{ll}
W(v,E), & \mbox{if } v \in E^+, \\
-W(v,E), & \mbox{if } v \in E^-, \\
0, & \mbox{otherwise.}
\end{array}\right.
\]

We observe that the multi-signed set model $(S,\mathcal{R})$ of a chemical reaction system defines a weighted hyperdigraph with incidence matrix $N$.  The construction proceeds as follows.  Let the chemical species of the chemical reaction system be denoted by the finite vertex set $S$.  We next obtain the hyperedge set $\mathcal{E}$ by defining a mapping from the collection of reactions $\mathcal{R}$ to $\{-1,0,+1\}^S$ given by $\mathcal{E}(r) = \sgn(N(\cdot,r))$ for $r \in \mathcal{R}$ so that $\mathcal{E}(\mathcal{R}) \subseteq \{-1,0,+1\}^S$.  Hence, each $r$ in $\mathcal{R}$ defines a directed hyperedge from which a definition for the weight function $W \in \mathbb{R}_+^{S \times \mathcal{R}}$ is given by setting 
\[
W(s,r) = |N(s,r)|.
\]
From this realization is follows that the incidence matrix $M \in \mathbb{R}^{S \times \mathcal{R}}$ is defined by setting
\[
M(s,r) = \mathcal{E}(r)(s)W(s,r) = \sgn(N(s,r))|N(s,r)| = N(s,r).
\]
Hence, the incidence matrix $N$ is shown to be sufficient to ``model'' the weighted hyperdigraph representation of the chemical reaction system defined by multi-signed sets on the set of chemical species.  In this approach, the directed hyperedges correspond to the set of reactions with the appropriate associated weights, which denote the change in molecularities that are induced by the given reaction.

Let $\mathcal{K} \in \mathbb{R}_{+}^{\mathcal{R}}$ be the kinetic rate constants for the set of chemical reactions $\mathcal{R}$.  Also let $X \in \mathbb{R}_+^S$ be the concentration of the chemical species $S$.  Now define a function $p \colon \mathcal{R} \to \mathbb{R}$ given by setting
\[
p(r) = \prod_{s \in S}X(s)^{\rho(r)(s)},
\]
for every $r$ in $\mathcal{R}$.  The function $p(\cdot)$ is defined to be the chemical potential \cite{weiss96}.  Next we define a function $J \colon \mathcal{R} \to \mathbb{R}$ by setting
\[
J(r) = \mathcal{K}(r)p(r),
\]
for all $r$ in $\mathcal{R}$ to define the flux.  With these definitions we next define the linear operator $\cdot$ obtained by setting
\[
\dot{X} = NJ,
\]
which corresponds to the differential equation derived from the dynamics of the chemical reaction network that was obtained from the constructed hyperdigraph model.

\begin{figure}[htbp]
\begin{center}
\subfigure[Enzymatic chemical reaction.]{
	\includegraphics[scale=.4]{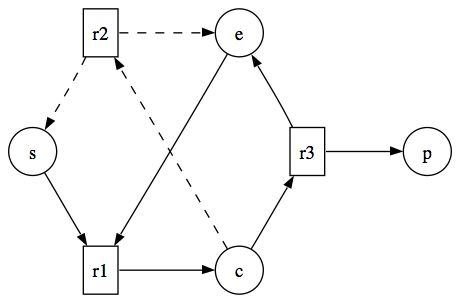}
	\label{fig:enzymeReaction}
	}
\hfill
\subfigure[A hypothetical example of a hyperdigraph.]{
	\includegraphics[scale=.4]{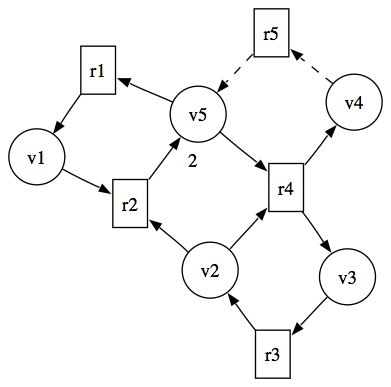}
	\label{fig:bailey}
	}
\caption{Two examples of hyperdigraphs.  The solid arrows indicate a hyperspanning tree.}
\end{center}
\end{figure}

The Michaelis-Menten ezymatic reaction \cite{murray02} will be considered for a small demonstration of the hyperdigraph model of chemical reaction networks.  The enzymatic reaction is kinetically represented as
\[
s + e \rlRev{r_1}{r_2} c \rTo{r_3} p + e,
\]
where reaction $r_1$ is complex formation of the substrate $s$ with the enzyme $e$ to produce the complex $c$; reaction $r_2$ is the disassociation of complex $c$ into the substrate and enzyme, which is also considered to be the reverse reaction of $r_1$; and reaction $r_3$ is the molecular reaction wherein the enzyme $e$ disassociates from the complex $c$ producing the product $p$.  The chemical species set is $S = \{s,e,c,p\}$ over which the collection of chemical complexes is defined as $\mathcal{C} = \{(1,1,0,0),(0,0,1,0),(0,1,0,1)\}$.  We further define the relation $\mathcal{R}$ representing the three reactions by $(1,1,0,0)\mathcal{R}(0,0,1,0)$, $(0,0,1,0)\mathcal{R}(1,1,0,0)$, and $(0,0,1,0)\mathcal{R}(0,1,0,1)$.  Finally, the incidence matrix is given by setting
\[
\begin{array}{rc}
& \begin{array}{ccc} r_1 & r_2 & r_3  \end{array} \\
N \quad = \quad \begin{array}{r}
s \\
e \\
c \\
p
\end{array} &
\left(\begin{array}{ccc}
-1 & 1 &  0  \\
-1 & 1 &  1  \\
1 & -1 & -1  \\
0 &  0 &  1
\end{array}\right).
\end{array}
\]
The hyperdigraph for this chemical reaction network is presented in Figure \ref{fig:enzymeReaction}, the vertices (chemical species) are denoted by circles and the directed hyperedges (chemical reactions) are denoted by rectangles.  The arrows denote the directionality of the hyperedge with respect each of its incident vertices.  Thus, the hyperedge labeled $r1$ represents transport of molecular mass from vertices $s$ and $e$ to vertex $c$.

We now construct an oriented matroid, which we refer to as a hypergraphic oriented matroid, on the hyperdigraph and show that the hypergraphic oriented matroid is isomorphic to the steady state solutions of the chemical reaction network.  It is worth noting that the hypergraphic oriented matroid corresponding to the steady state solutions of the dynamical system are equivalent to the $S$-invariants of the Petri net chemical reaction model \cite{goss98}.  This equivalence is due to the fact that hyperdigraphs are isomorphic to Petri nets \cite{oliveira03}.

\section{Hypergraphic Oriented Matroids}
Let $S$ be a finite set and consider the multisets $x \in \mathbb{Z}^S$ defined on $S$, where $x(s)$ is the multiplicity of the element $s$ in $S$.  Denote $\zero$ to be the multiset that maps every element $s$ to $0$.  We observe that allowing negative multiplicities is a generalization of multisets \cite{reisig85}, which are traditionally defined by mappings to the natural numbers \cite{blizard89, calude01}.  We will define the operation of addition on $\mathbb{Z}^S$ where, for arbitrary multisets $x$ and $y$ in $\mathbb{Z}^S$ and $s$ in $S$ we have $(x + y)(s) = x(s) + y(s)$.  Further, we define the unary operation of negation as $(-x)(s) = -x(s)$.  It is now clear that $(\mathbb{Z}^S, +)$ defines an additive abelian group.  Further we may define scalar multiplication as the function $\mathbb{Z}\times \mathbb{Z}^S \mapsto \mathbb{Z}^S$, given by $(\alpha x)(s) = \alpha x(s)$ for arbitrary $\alpha$ in $\mathbb{Z}$, $x$ in $\mathbb{Z}^S$, and $s$ in $S$.  We observe that for $\alpha$ and $\beta$ in $\mathbb{Z}$, and $x$ and $y$ in $\mathbb{Z}^S$, we have the following distributive properties $(\alpha + \beta)x = \alpha x + \beta x$ and $\alpha (x + y) = \alpha x + \alpha y$.  We may therefore consider $\mathbb{Z}^S$ as a module which is referred to as the $\mathbb{Z}$-module \cite{maclane99,roman05}.  

Recall the multiset model of the chemical reaction system defined previously, wherein the set $S$ is the set of chemical species and the matrix $N$ represented a finite collection of multisets $N(\cdot, \mathcal{R})$ on $S$.  Given that the set $\mathcal{R}$ corresponds to the hyperedges of the hyperdigraph, we have that the two dual oriented matroids correspond, respectively, to the linearly independent and linearly dependent subsets of $\mathcal{R}$.  The minimal linearly dependent subsets of $\mathcal{R}$ will correspond to hypercycles of the hyperdigraph.  Similarly, the maximal linearly independent subsets of $\mathcal{R}$ will correspond to the hyperspanning forests of the hyperdigraph.

We will now proceed to construct the hypergraphic oriented matroid on the $\mathbb{Z}$-module $M \subseteq \mathbb{Z}^\mathcal{R}$.  An operator $\mathrm{cl}$ is on the $\mathbb{Z}$-module $M$ as the function from $2^M$ into $2^M$ given by, for any $X \subseteq M$,
\[
\mathrm{cl}(X) = \left\{ m \in M ~\left|~ bm = \sum_{x \in X}\alpha(x)x, \mbox{ for } \alpha \in \mathbb{Z}^X, ~ 0 \neq b \in \mathbb{Z}\right.\right\}.
\]
It was shown by \cite{danielsson03} that the function $\mathrm{cl}$ is a closure operator  and so satisfies the following properties \cite{oxley05},
\begin{itemize}
\item[(CL1)] if $X \subseteq M$, then $X \subseteq \mathrm{cl}(X)$;
\item[(CL2)] if $X \subseteq Y \subseteq M$, then $\mathrm{cl}(X) \subseteq \mathrm{cl}(Y)$;
\item[(CL3)] if $X \subseteq M$, then $\mathrm{cl}(\mathrm{cl}(X)) = \mathrm{cl}(X)$;
\item[(CL4)] if $X \subseteq M$, $x \in M$, and $y \in \mathrm{cl}(X\cup x)\setminus \mathrm{cl}(X)$, then $x \in \mathrm{cl}(X \cup y)$.
\end{itemize}
We observe that for $X$ subset of $M$, $\mathrm{cl}(X)$ is a submodule of the $\mathbb{Z}$-module $M$ \cite{danielsson03}.

For the construction of the matroid bases it is important to define an irreducible multiset.  To do so, we will define a function $g \colon M \to \mathbb{Z}$ given by setting, for any $x$ in $M$,
\[
g(x) = \sum_{r \in \mathcal{R}}\alpha(r)x(r),
\]
where $\alpha$ in $\mathbb{Z}^\mathcal{R}$ is chosen such that $g(x) = \gcd\{x(r) ~|~ r \in \mathcal{R}\}$.  We then define the mapping $\mathring{~} \colon M \to M$, which will be referred to as the ``reducing map" by setting
\[
\mathring{x}(r) = \frac{1}{g(x)}x(r),
\]
for every $r$ in $\mathcal{R}$.  A multiset $x$ is irreducible if $x = \mathring{x}$.  The intention of the reducing map is to obtain an analog of mapping a vector in a vector field to its corresponding unitary vector.  Define $\zero$ to be irreducible.  For a subset of $X$ of $M$, we define $\mathring{X} = \{\mathring{x} \in M ~| ~ x \in X\}$.  

We now present an algorithm similar to Fourier-Motzkin \cite{pfeiffer99} to construct bases for the dual hypergraphic oriented matriods $M$ and $M^*$ of $\mathbb{Z}^\mathcal{R}$.  

\begin{alg}\label{alg:fourierMotzkin}
{\bf Basis:} Let $S$ be a finite set and let $\mathcal{R}$ be a finite subset of $\mathbb{Z}^S$, which will be represented by the integral matrix $N \in \mathbb{Z}^{S \times \mathcal{R}}$.  We will construct the mappings $F \in \mathbb{Z}^{\mathcal{R} \times (S \disjointunion \mathcal{R})}$ and $F^* \in \mathbb{Z}^{S \times (\mathcal{R} \disjointunion S)}$ by
\[
\begin{array}{rc}
& \begin{array}{ccc} S & & \mathcal{R}  \end{array} \\
F \quad =  & \mathcal{R}
\left(\begin{array}{cc}
N^T &  \mathrm{\bf Id}(\mathcal{R})
\end{array}\right),
\end{array} \quad 
\begin{array}{rc}
& \begin{array}{lcr} \mathcal{R} & & S  \end{array} \\
F^* \quad =  & S
\left(\begin{array}{cc}
N &  \mathrm{\bf Id}(S)
\end{array}\right),
\end{array}
\]
where $\mathrm{\bf Id}(\mathcal{R})$ and $\mathrm{\bf Id}(S)$ are identity matrices. We will then perform Gaussian elimination on $F$ and $F^*$ where the row elimination process of using $F(r_i,\cdot)$ to eliminate $F(r_j,\cdot)$ is
\begin{eqnarray*}
c & := & \lcm\{F(r_i,s), F(r_j,s)\} \\
a & := & c/F(r_i,s) \\
b & := & c/F(r_j,s) \\
F(r_j,\cdot) & := & bF(r_j,\cdot) - aF(r_i,\cdot)
\end{eqnarray*}
After the completion of Gaussian elimination on the matrices $F$ and $F^*$, each is partitioned into submatrices as
\[
\begin{array}{rc}
& \begin{array}{ccc} S & \qquad & \mathcal{R}  \end{array} \\
F \quad =  & \mathcal{R}
\left(\begin{array}{cc}
\im(N) &  \mathrm{Dom}/\ker(N) \\
\zero & \ker(N) 
\end{array}\right),
\end{array}
\]
\[
\begin{array}{rc}
& \begin{array}{ccc} \mathcal{R} & \qquad & S  \end{array} \\
F^* \quad =  & S
\left(\begin{array}{cc}
\im(N^T) &  \mathrm{Dom}/\ker(N^T) \\
\zero & \ker(N^T) 
\end{array}\right).
\end{array}
\]
The respective basis sets for $M$ and $M^*$ are then defined by the submatrix of $F$ labeled $\ker(N)$ and the submatrix of $F^*$ labeled $\im(N^T)$.  That is,
\[
B = \{F(r,\mathcal{R}) ~|~ F(r,S) = 0, ~r \in \mathcal{R}\},
\] 
and
\[
B^* = \{F^*(s,\mathcal{R}) ~|~ F^*(s,\mathcal{R}) \neq 0, ~s \in S\},
\] 
We thus obtain the sets $M = \mathring{(\mathrm{cl}(B))}$ and $M^* = \mathring{(\mathrm{cl}(B^*))}$, which correspond to the desired hypergraph oriented matroids.
\end{alg}

We remark that even though $B$ may not be a basis set for the oriented matroid, it is sufficient to identify a fundamental set of independent hypercycles as well as the matroid rank dimension.  As such, future considerations may be to determine whether $B$ is in fact a basis for $M$ and also check whether $\mathring{X} \subseteq \mathrm{cl}(X)$.

\section{Hypercycles and Hyperspanning Trees}
The dual matroids $M$ and $M^*$ correspond respectively to the cycle space and co-cycle space of the chemical reaction network represented by the incidence matrix $N$.  The elements of $M$ are hypercycles, the basis elements of which are minimal linearly dependent subsets of the reactions that define the hyperedges of the network.  We use the term ``hypercycle" to make a distinction between the elements of $M$ and cycles in hypergraphs as defined by Berge \cite{berge73}.  

In the general theory of hypergraphs \cite{berge73}, a chain of length $q$ in a hypergraph is a sequence $v_1 \rTo{r_1} v_2 \rTo{r_2} \dots \rTo{r_q} v_{q+1}$ such that
\begin{itemize}
\item[(C1)] $v_1,v_2,\dots,v_{q}$ are all distinct vertices of $\mathcal{H}$;
\item[(C2)] $r_1,r_2,\dots,r_q$ are all distinct hyperedges of $\mathcal{H}$;
\item[(C3)] $v_k, v_{k+1} \in \supp(r_k)$ with $\rho(r_k)(v_k)\rho(r_k)(v_{k+1}) = 0 = \pi(r_k)(v_k)\pi(r_k)(v_{k+1})$ for $k=1,\dots,q$.
\end{itemize}
If $q > 1$ and $v_{q+1} = v_1$, then this chain is called a cycle of length $q$ \cite{berge73}.  This general definition of cycle corresponds to the cycles which exist in the bipartite representation of the hypergraph, and therefore does not incorporate the nonlinear dependencies that define the hyperedges in chemical reaction networks.  For this reason we will use the term ``closed loop" in leu of ``cycle" for this class of hyperdigraphs and define the more restrictive hypercycle.

A hypercycle of length $q$ is a multiset $y \in \mathbb{Z}^\mathcal{R}$ such that $N(s,\cdot) \perp y$ for every $s$ in $S$ with $|\supp(y)| = q$.  Observe that $y(r)$ is the resulting mass flux caused by reaction $r$ in $\mathcal{R}$.  Further, there is not an order restriction on the set of hyperedges that construct the hypercycle.  Since $N(s,\cdot) \cdot y = 0$ is the net change in mass for chemical species $s$ resulting from $y$, the hypercycles, and thus the matroid $M$, represent steady state flows in the network.

A subhyperdigraph of a hyperdigraph $\mathcal{H} = (V,\mathcal{R})$ is a hyperdigraph $\mathcal{H}' = (V',\mathcal{R}')$ such that $V'$ and $\mathcal{R}'$ are subsets of $V$ and $\mathcal{R}$, respectively, and each hyperedge in $\mathcal{R}'$ is a multiset on $V$ with $N(V\setminus V', r) \mapsto 0$ for every $r \in \mathcal{R}'$.  A hyperspanning forest in a hyperdigraph $\mathcal{H}$ is a maximal subhyperdigraph $\mathcal{H}' = (V',\mathcal{R}')$ such that $\mathcal{H}'$ does not contain any hypercycles and $V = V'$.  

\begin{cor}\label{cor:cyclomaticnumber}
The hypercyclomatic number for a hypergraph $\mathcal{H} = (V,\mathcal{R})$ is 
\[
c = \rank(M) = \nullity(N) = |\mathcal{R}| - \rank(N^T),
\]
where $M$ and $M^*$ are the matroids for $\ker(N)$ and $\im(N^T)$, respectively.  
\end{cor}

In general there may exist closed loops in hyperspanning trees as indicated in Figure \ref{fig:bailey}, where the hyperspanning tree is depicted by the solid arrows.  There are three independent closed loops in the hyperspanning tree, namely 
\[
\begin{array}{c}
v_1 \rTo{r_2} v_5 \rTo{r_1} v_1, \\
v_5 \rTo{r_4} v_3 \rTo{r_3} v_2 \rTo{r_2} v_5, \\
v_2 \rTo{r_4} v_3 \rTo{r_3} v_2.
\end{array}
\]
The only hypercycle in this example is $y := \langle 0,0,1,1,1\rangle$.  The incidence matrix $N$ is 
\[
\begin{array}{rc}
& \begin{array}{ccccc} r_1 & r_2 & r_3 & r_4 & r_5  \end{array} \\
N \quad = \quad \begin{array}{r}
v_1 \\
v_2 \\
v_3 \\
v_4 \\
v_5
\end{array} &
\left(\begin{array}{ccccc}
1 & -1 &  0 & 0  & 0 \\
0 & -1 &  1 & -1 & 0 \\
0 & 0 & -1  & 1  & 0 \\
0 &  0 &  0 & 1  & -1 \\
-1 & 1 & 0 & -1 &  1
\end{array}\right).
\end{array}
\]

\section{MAP Kinase Cascade}
The MAP Kinase cascade is composed of six coupled enzymatic reactions and a complex formation \cite{bhalla99}, which are listed in the following tableau, in which in the interest of space we have abbreviated the enzymatic reactions $S + E \rlRev{}{} C \rTo E + P$ as $S \rTo{E} P$.  
\[
\begin{array}{l|l}
Raf \rlRev{PKC}{PP2-A} Raf^*	& Raf^* \rlRev{MAPK^*}{PP2-A} Raf^{**} \\
MAPKK \rlRev{GTP.Ras.Raf^*}{PP2-A} MAPKK^* & MAPKK^* \rlRev{GTP.Ras.Raf^*}{PP2-A} MAPKK^{**} \\
MAPK \rlRev{MAPKK^{**}}{MKP1} MAPK_{tyr^*} & MAPK_{tyr^*} \rlRev{MAPKK^{**}}{MKP1} MAPK_{tyr^*} \\
Raf^* + GTP.Ras \rlRev{}{} GTP.Ras.Raf^*	&
\end{array}
\]
\begin{figure}[htbp]
\begin{center}
\includegraphics[scale=.2]{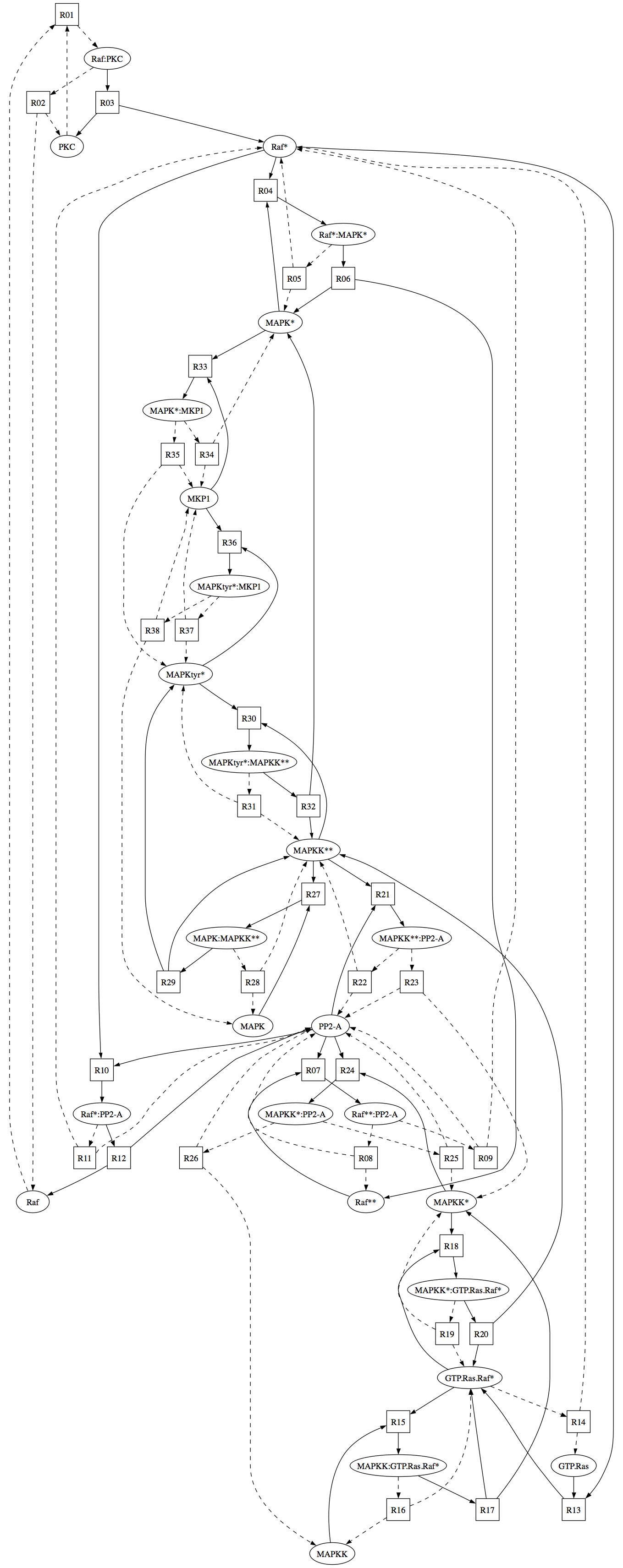}
\caption{Hyperdigraph of the MAP Kinase Cascade.  The solid arrows indicate one of the many possible hyperspanning trees.}
\label{fig:mapk}
\end{center}
\end{figure}
The hyperdigraph representation of the MAP Kinase cascade is presented in Figure \ref{fig:mapk}, with the solid lines indicating one of many possible the hyperspanning trees.  The rank of the hypercycle, hypergraphic oriented matroid is $19$, where there are three irreducible hypercycles for each of the six coupled enzymatic reactions i) $S \rlRev{E_1}{E_2} P$, ii) $S + E_1 \rlRev{}{} S\colon E_1$, and iii) $P + E_2 \rlRev{}{} P\colon E_2$; and one irreducible hypercycle for the reaction $Raf^* + GTP.Ras \rlRev{}{} GTP.Ras.Raf^*$.  Recall that the hypercycles correspond to the steady state solutions of the dynamic system.  As such, any linear combination of the $19$ hypercycles is also a steady state solution.  The dual hypergraphic oriented matroid, which corresponds to the hyperspanning trees, also has rank $19$.  We observe that there are $8$ closed loops within the hyperspanning tree presented in Figure \ref{fig:mapk}.  Six of the eight closed loops are linearly independent and are of length $2$.  The remaining two closed loops are each of length $4$ and are linearly dependent on the six closed loops.  Each of the six linearly independent closed loops correspond to one-half of a coupled enzymatic reaction, wherein the closed loop explicitly denotes the conservation of the enzyme while producing the product from the substrate as in the general depiction $S + E \rTo{} S\colon E \rTo{} E + P$.  The two linearly dependent closed loops of length $4$ correspond to a conservation of the enzyme in a sequence of enzymatic reactions in which both reactions are catalyzed by the same enzyme as in the general depiction $S + E \rTo{} S\colon E \rTo{} E + S_0 \rTo{} S_0\colon E \rTo{} P + E$. 

Oliveira \emph{et al.} \cite{oliveira01,oliveira03,oliveira04} proposed and successfully implemented a method for rank ordering the chemical species and reactions within a chemical reaction system by counting the occurrences of each species (reaction) within the complete list of closed loops (Berge cycles) contained in the hyperdigraph of the chemical reaction network.  It is observed that rank ordering the chemical species based on the complete listing of closed loops is analogous to a rank ordering based on the centrality of the chemical species within the hyperdigraph \cite{estrada06}.  Estrada and Rodr\'{i}guez-Vel\'{a}quez \cite{estrada06} characterizes the centrality rank index based on the hypergraph adjacency matrix.  Using the hyperdigraph model proposed in this paper we can construct the hyperdigraph adjacency matrix for the chemical reaction network by setting $L = A^TB$, where the matrices $A$ and $B$ were defined previously.  

\begin{figure}[htbp]
\begin{center}
\includegraphics[scale=.4]{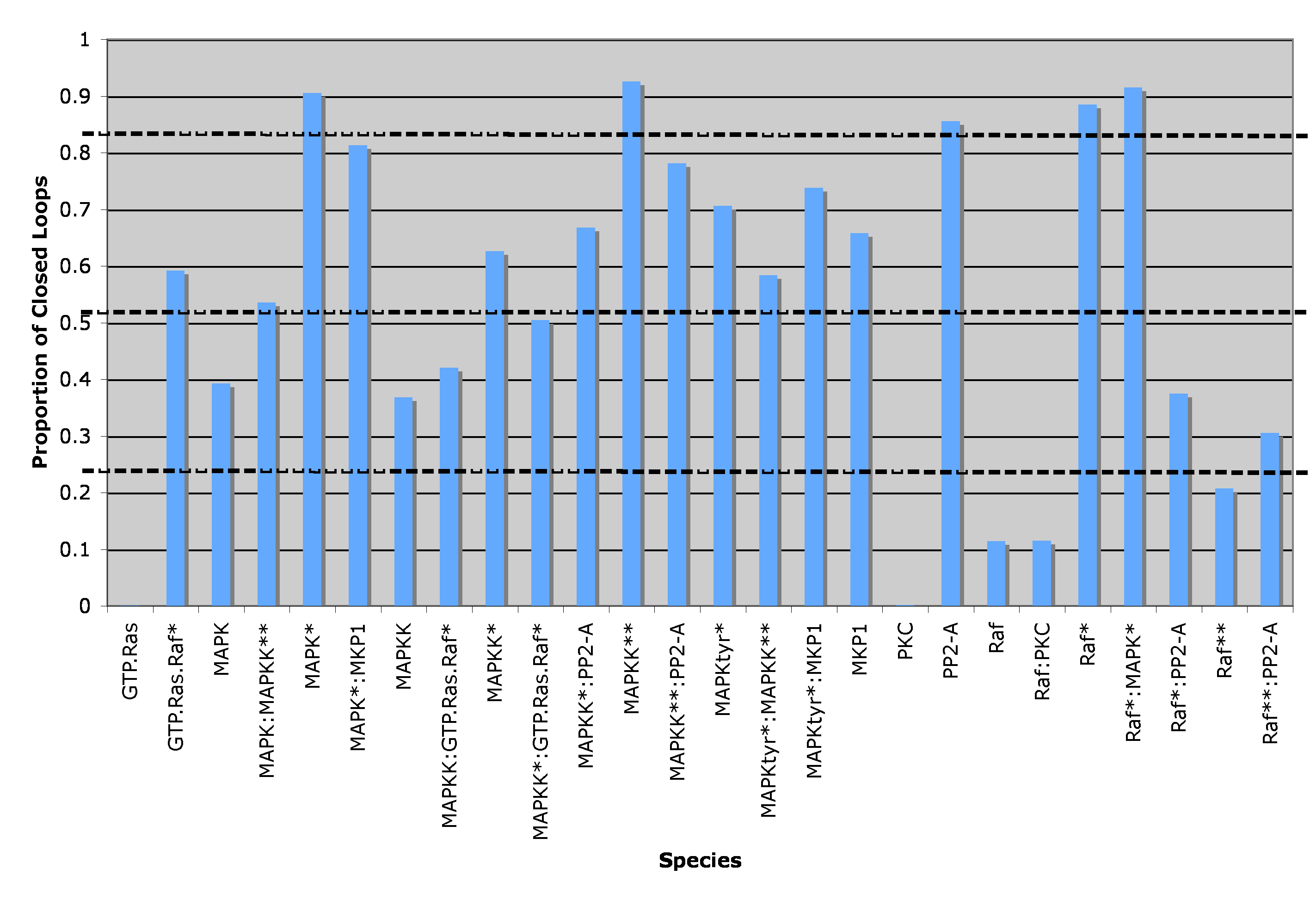}
\caption{Centrality ranking index of Map Kinase cascade protein species.  The dashed, horizontal lines represent the mean + standard deviation, mean, and mean - standard deviation, respectively.}
\label{fig:speciesCentrality}
\end{center}
\end{figure}

There are $1456$ irreducible closed loops within the MAP Kinase cascade hyperdigraph shown in Figure \ref{fig:mapk}.  Figure \ref{fig:speciesCentrality} shows the proportion of the $1456$ irreducible closed loops incident with each of the chemical species.  For the chemical species the mean proportion is $0.538$ and the standard deviation is $0.289$.  Oliveira \emph{et al.} futher demonstrated that chemical species which occur most and least often are critical to the functionality of the chemical reaction network.  Those species that occur most often within the closed loops, that is being most central, may be considered critical to the network as a ``pinch point" or being utilized often in the signaling pathway.  Those species that occur least often within the closed loops, that is being least central, may be considered critical to the network as an initiator or trigger for the cell-signaling process.  

For this analysis the designations of high and low centrality will respectively be considered by the values of the mean + standard deviation ($0.827$) and the mean - standard deviation ($0.249$).  The most central chemical species are $MAPKK^{**}$ ($0.926$), $Raf^*\colon MAPK^*$ ($0.915$), $MAPK^*$ ($0.905$), $Raf^*$ ($0.885$), and $PP2-A$ ($0.856$).  The least central chemical species are $GTP.Ras$ ($6.87\times 10^{-4}$), $PKC$ ($1.37\times 10^{-3}$), $Raf$ ($0.114$), $Raf\colon PKC$ ($0.115$), and $Raf^{**}$ ($0.207$).

\section{Discussion}
We have effectively demonstrated the construction of a multiset model for chemical reaction networks, which canonically generates a hyperdigraph.  The hyperdigraph model was utilized in the qualitative steady state analysis of the Map Kinase cascade, wherein it was shown that each of the six coupled enzymatic reactions contribute three independent hypercycles with an additional hypercycle corresponding to the complex formation $Raf^* + GTP.Ras \rlRev{}{} GTP.Ras.Raf^*$.  The hypercycles are invariants in the dynamical system of the chemical reaction network. 

We further addressed the question of a centrality rank index which rank orders the vertices in a hyperdigraph model of chemical reaction networks.  This centrality rank index highlighted ten protiens of interest in the Map Kinase cascade, five of which are highly utilized in the Map Kinase cascade while the remaining five are most likely signaling pathway regulators.

The authors would like to Professor Hong Qian, Dr. Janet B. Jones-Oliveira, Professor Stefan E. Schmidt, and Dr. Thorsten Pfeiffer for their wonderful insights and discussions.  The figures in this paper where generated by Graphviz version 1.13(16) a product of AT\&T and Pixelglow Software (http://www.pixelglow.com/graphviz).

\bibliography{bibliography}
\end{document}